\documentclass[11pt]{amsart}
\usepackage{amsmath}
\usepackage{amssymb}
\usepackage{amsfonts}
\usepackage{graphicx}

\DeclareSymbolFont{AMSb}{U}{msb}{m}{n}
\DeclareMathSymbol{\N}{\mathbin}{AMSb}{"4E}
\DeclareMathSymbol{\Z}{\mathbin}{AMSb}{"5A}
\DeclareMathSymbol{\R}{\mathbin}{AMSb}{"52}
\DeclareMathSymbol{\Q}{\mathbin}{AMSb}{"51}
\DeclareMathSymbol{\I}{\mathbin}{AMSb}{"49}
\DeclareMathSymbol{\C}{\mathbin}{AMSb}{"43}

\theoremstyle{definition}

\theoremstyle{corollary}

\theoremstyle{example}

\theoremstyle{note}

\theoremstyle{notation}

\numberwithin{equation}{section}
\begin{document}
\title{MONOCHROMATIC SUMS EQUAL TO PRODUCTS NEAR ZERO}
\author{Sourav Kanti Patra and Md Moid Shaikh}

\address{Department of Mathematics, Ramakrishna Mission Vidyamandira,
Belur Math, Howrah-711202, West Bengal, India}
\email{souravkantipatra@gmail.com}

\address{Department of Mathematics, Maharaja Manindra Chandra College,
20, Ramkanto Bose Street, Kolkata-700 003, West Bengal, India}
\email{mdmoidshaikh@gmail.com}

%\author{}
%\address{}
%\email{}
%\thanks{}

\begin {abstract} 
Hindman proved that, whenever the set $\mathbb {N}$ of naturals is finitely 
colored, there must exist non-constant monochromatic solution of the 
equation $a+b=cd$. In this paper we extend this result for dense subsemigroups of 
$((0, \infty), +)$ to near zero.

\end{abstract}

\maketitle

\section {introduction}
In Ref. \cite{csikvari}, P. Csikv\'{a}ri, A. S\'{a}rk\"{o}zy, and K. Gyarmati asked whether
whenever the set $\mathbb {N}$ of naturals is finitely 
colored, there must exist monochromatic $a,b,c$ and $d$ with $a\neq b$ such that $a+b=cd$.
In Ref. \cite{Hindman1}, Hindman answered this question affirmatively
by showing in addition that one can demand that $a,b, c, d$ are all distinct and the color of
$a+b$ is that same as that of $a,b,c$ and $d$. In fact he proved considerably more stronger 
result using the algebraic structure of $\beta \mathbb {N}$, the Stone-$\breve{C}$ech compactification
of $\mathbb {N}$.

Let $(S,\cdot)$ be an infinite discrete semigroup. Now the points
of $\beta S$ are taken to be the ultrafilters on $S$, the
principal ultrafilters being identified with the points of $S$.
Given $A\subseteq S$, let us set $\bar{A}=\{p\in\beta S:A\in p\}$.
Then the set $\{\bar{A}:A\subseteq S\}$ will become a basis for a
topology on $\beta S$. The operation $\cdot$ on $S$ can be
extended to the Stone-$\breve{C}$ech compactification $\beta S$ of
$S$ so that $(\beta S,\cdot)$ is a compact right topological
semigroup (meaning that for any $p\in \beta S$, the function
$\rho_p:\beta S \rightarrow \beta S$ defined by $\rho_p(q)=q\cdot
p$ is continuous) with $S$ contained in its topological center
(meaning that for any $x\in S$, the function $\lambda_x:\beta S
\rightarrow \beta S$ defined by $\lambda_x(q)=x\cdot q$ is
continuous). Given $p,q\in\beta S$ and $A\subseteq S$, $A\in p\cdot
q$ if and only if $\{x\in
S:x^{-1}A\in q\}\in p$, where $x^{-1}A=\{y\in S:x\cdot y\in A\}$.\\

A nonempty subset $I$ of a semigroup $(T,\cdot)$ is called a left
ideal of $T$ if $T\cdot I\subseteq I$, a right ideal if
$I.T\subseteq I$, and a two-sided ideal (or simply an ideal) if it
is both a left and a right ideal. A minimal left ideal is a left
ideal that does not contain any proper left ideal. Similarly, we
can define minimal right ideal and smallest ideal. Any compact
Hausdorff right topological semigroup $(T,\cdot)$ has the unique
smallest two-sided ideal
$$\begin{array}{ccc}
K(T) & = & \bigcup\{L:L \text{ is a minimal left ideal of } T\} \\
& = & \,\,\,\,\,\bigcup\{R:R \text{ is a minimal right ideal of } T\}\\
\end{array}$$

Given a minimal left ideal $L$ and a minimal right ideal $R$ of
$T$, $L\cap R$ is a group, and in particular $K(T)$ contains an
idempotent. An idempotent that belongs to $K(T)$ is called a
minimal idempotent.\\

In Ref. \cite{Hindman3}, authors introduced different notations of largeness for 
arbitrary semigroup and studied their combinatorial properties. From now we use $\mathcal{P}_{f}(X)$ 
to denote the set of all finite subsets of a set $X$.

\textbf{Definition 1.1} \cite[Definition 3.1]{Hindman3}: Let $(S,\cdot)$ be a semigroup.

(a)  A set $A\subseteq S$ is syndetic if only if there exists $G \in \mathcal{P}_{f}(S)$ with

$S\subseteq \bigcup_{t\in G}  t^{-1}A$.

(b)  A set $A\subseteq S$ is  piecewise syndetic if only if there exists $G \in \mathcal{P}_{f}(S)$ such
that $\{y^{-1}(\bigcup_{t\in G}  t^{-1}A) :y\in S\}$ has the finite intersection property.

(c) A family $\mathcal{A}\subseteq \mathcal{P}(S)$ is collectionwise piecewise syndetic 
if and only if there exists a function $G: \mathcal{P}_{f}(\mathcal{A})\longrightarrow\mathcal{P}_{f}(S)$
such that

$\{y^{-1}(G(\mathcal{F}))^{-1}(\cap\mathcal{F}): y\in S$  and $\mathcal{F}\in \mathcal{P}_{f}(\mathcal{A}) \}$
has the finite intersection property.

\textbf{Theorem 1.2}: Let $(S,\cdot)$ be a semigroup and let $\mathcal{A}\subseteq \mathcal{P}(S)$.
Then there exists $p\in K(\beta S)$ with $\mathcal{A}\subseteq p$ if and only if $\mathcal{A}$ is 
collectionwise piecewise syndetic. In particular, given $A\subseteq S$, $K(\beta S)\cap cl A\neq\O{}$ if and only 
if $A$ is piecewise syndetic.

\textbf{Proof}: See the proof \cite[Theorem 3.2]{Hindman3}.\\

Furstenberg introduced a class of large sets in terms of notion from topological dynamics.
There is a nice characterization of such sets in terms of algebraic structure of $\beta\mathbb{N}$.

We now recall the following definitions \cite [Definitions 4.42 and 15.3]{Hindman4}

\textbf{Definition 1.3}: Let $(S,\cdot)$ be a semigroup and let $A\subseteq S$.

(a) $A$ is called Central in $(S,\cdot)$ if there is some idempotent $p\in K(\beta S,\cdot) $
such that $A\in p$.

(b) $A$ is called Central* in $(S,\cdot)$ if $A\cap B \neq\O{}$ for every central set $B$ in
$(S,\cdot)$.

In Ref. \cite{Bergelson1}, authors investigated the interplay between additive and multiplicative
largeness. We need the following definition \cite[Definition 1]{Hindman1} to 
state such results. 

\textbf{Definition 1.4}: Let $\langle{x_n}\rangle_{n=1}^{\infty}$ be an infinite sequence of
positive real numbers, let $m\in\mathbb{N}$ 
and let $\langle{y_n}\rangle_{n=1}^{m}$be a finite sequence of positive real numbers. Then\\

(a)   $FS(\langle{x_n}\rangle_{n=1}^{\infty})=\{\sum_{n\in F}x_n:F\in\mathcal{P}_f(\mathbb N)\}$ and\\
    
      $FP(\langle{x_n}\rangle_{n=1}^{\infty})=\{\prod_{n\in F}x_n:F\in\mathcal{P}_f(\mathbb N)\}$.\\

(b)  $FS(\langle{y_n}\rangle_{n=1}^{m})=\{\sum_{n\in F}y_n:\O{}\neq F\subseteq\{1,2,.....,m\}\}$ and\\ 

     $FP(\langle{y_n}\rangle_{n=1}^{m})=\{\prod_{n\in F}y_n:\O{}\neq F\subseteq\{1,2,.....,m\}\}$.\\
     
(c)  The sequence $\langle{x_n}\rangle_{n=1}^{\infty}$ has distinct finite sums if and only if whenever

$F, G\in\mathcal{P}_f(\mathbb N)$ and $F\neq G$ , one has $\sum_{n\in F}x_n\neq \sum_{n\in G}x_n$. 
The analogous definition applies to $\langle{y_n}\rangle_{n=1}^{m}$.

(d) The sequence $\langle{x_n}\rangle_{n=1}^{\infty}$ has distinct finite products if and only if
whenever $F, G\in\mathcal{P}_f(\mathbb N)$ and $F\neq G$ , one has $\prod_{n\in F}x_n\neq \prod_{n\in G}x_n$.
The analogous definition applies to $\langle{y_n}\rangle_{n=1}^{m}$.

(e) The sequence $\langle{x_n}\rangle_{n=1}^{\infty}$ is strongly increasing if and only if for each 
    $n\in \mathbb{N}$, $\sum_{t=1}^nx_{t}< x_{n+1}$.
    
 Clearly if  $\langle{x_n}\rangle_{n=1}^{\infty}$ is strongly increasing, then it has distinct finite sums. 
 We also need to recall the following definition from Ref. \cite{Bergelson1}.
 
\textbf{Definition 1.5}: Let $A\subseteq \mathbb N$. Then $A$ is said to be an $IP_0$-set in $(\mathbb{N}, +)$ 
if and only if for each $m\in\mathbb N$, there exists a finite sequence 
$\langle{y_n}\rangle_{n=1}^{m}$ such that $FS(\langle{y_n}\rangle_{n=1}^{m})\subseteq A$.

\textbf{Theorem 1.6}: For all $A\subseteq \mathbb N$, if $A$ is syndetic in $(\mathbb{N},\cdot)$ then $A$ is
central in $(\mathbb{N},+)$.\\
\textbf{Proof}: See the proof \cite[Lemma 5.11]{Bergelson2}.

\textbf{Theorem 1.7}: For all $A\subseteq \mathbb N$, if $A$ is central in $(\mathbb{N},\cdot)$ then 
 $A$ is an $IP_0$-set in $(\mathbb{N},+)$.\\ 
\textbf{Proof}: See the proof \cite[Theorem 3.5]{Bergelson3}.\\

   In the present paper we want to extend the Theorem $5$ of Ref. \cite{Hindman1} in section 3 for a dense
   subsemigroup of $((0, \infty), +)$. Let $S$ be a dense subsemigroup of $((0, \infty), +)$, then one can 
   define \\
   $ 0^+(S)= \{p\in \beta S_d:$ for all $(\epsilon>0)((0, \epsilon), \in p)\}$,
   ($S_d$ is the set $S$ with the discrete topology).\\
   
   \textbf{Theorem 1.8}: Let $S$ be a dense subsemigroup of $((0, \infty), +)$ such that $S\cap(0,1)$
is a subsemigroup of $((0, 1), \cdot)$ and assume that for each $y\in S\cap(0,1)$ and each $x\in S$, $x/y$ and
$yx\in S$. Then $ 0^+(S)$ is a two-sided ideal of $\beta(S\cap(0,1),\cdot)$.\\
\textbf{Proof}: See the proof\cite[Lemma 13.29.(f)]{Hindman4} 

 \section {Additive and Multiplicative Largeness Near Zero}
 Hindman and Leader introduced different notions of large sets near zero. In this section we study 
 the interplay between additive and multiplicative large sets near zero.
 We now recall some notions from Ref. \cite{Hindman2}. \\
 \textbf{Definition 2.1}: Let $S$ be a dense subsemigroup of $((0, \infty), +)$.\\
 (a) \cite[Definition 3.2]{Hindman2}: A set $A\subseteq S$ is syndetic near zero if only if 
 for every $\epsilon>0$ there exist some 
 $F \in \mathcal{P}_{f}((0, \epsilon)\cap S)$ and some $\delta>0$ such that
 $S\cap (0, \delta)\subseteq \bigcup_{t\in F} (-t+A)$.\\
 (b) \cite[Definition 3.4]{Hindman2}: A subset $A$ of $S$ is piecewise syndetic near zero if and only if 
 there exist sequences $\langle{F_n}\rangle_{n=1}^{\infty}$ and $\langle{\delta_n}\rangle_{n=1}^{\infty}$ such that\\
 (1) for each $n\in\mathbb{N}$, $F_n\in\mathcal{P}_{f}((0, 1 /n)\cap S)$ and $\delta_n\in (0, 1 /n)$ and\\
 (2) for all $G\in\mathcal{P}_{f}(S)$ and all $\mu>0$ there is some $x\in(0, \mu)\cap S$ such that for all
    $n\in\mathbb{N}$, $(G\cap (0, \delta_{n}))+x\subseteq \bigcup_{t\in F_n} (-t+A)$.\\
    
\textbf{Theorem 2.2}  \cite[Theorem 3.5]{Hindman2}: Let $S$ be a dense subsemigroup of $((0, \infty), +)$
and let $A\subseteq S$. Then $K(0^+(S))\cap\bar{A} \neq\O{}$ if and only if  $A$ is piecewise syndetic near zero.\\

\textbf{Definition 2.3}: Let $S$ be a dense subsemigroup of $((0, \infty), +)$. A subset  $A$ of $ S$ is said to be an 
$IP$-set near zero if and only if for each $\epsilon>0$ there exists some sequence $\langle{x_n}\rangle_{n=1}^{\infty}$ in $S$
such that $\sum_{n=1}^{\infty}x_{n}$ converges and
$FS(\langle{x_n}\rangle_{n=1}^{\infty})\subseteq A\cap (0, \epsilon)$.\\

\textbf{Theorem 2.4}: Let $S$ be a dense subsemigroup of $((0, \infty), +)$ and let  $A\subseteq S$. Then $A$
is an $IP$-set near zero  if and only if there is some idempotent $p$ in $0^+(S)$ such that $A\in p$.\\
\textbf{Proof}: See the proof \cite[Theorem 3.1]{Hindman2}.\\

\textbf{Definition 2.5}: Let $S$ be a dense subsemigroup of $((0, \infty), +)$. 
Define $\varGamma_0(S)=$\{$p\in \beta S:$ if  $A\in p$ then $A$ is an $IP$-set near zero\}\\

\textbf{Lemma 2.6}: Let $S$ be a dense subsemigroup of $((0, \infty), +)$ such that $S\cap(0,1)$
is a subsemigroup of $((0, 1), \cdot)$ and assume that for each $y\in S\cap(0,1)$ and each $x\in S$, $x/y$ and
$yx\in S$. Then $\varGamma_0(S)$ is a left ideal of $\beta(S\cap(0,1),\cdot)$.\\

\textbf{Proof}: Note that $\varGamma_0(S)= cl\{p\in 0^+(S): p+p=p\}$, and therefore is non-empty.
Let $p\in \varGamma_0$, and $q\in \beta(S\cap(0,1),\cdot)$, and let $A\in q\cdot p$. 
Then $\{y\in S\cap(0,1):y^{-1}A\in p\}\in q$, so pick $y\in S\cap(0,1)$ such that $y^{-1}A\in p$.
Thus $y^{-1}A$ is an $IP$-set near zero. Take some sequence $\langle{x_n}\rangle_{n=1}^{\infty}$
in $S$ such that $FS(\langle{x_n}\rangle_{n=1}^{\infty})\subseteq (y^{-1} A)\cap (0, \epsilon)$ and 
$\sum_{n=1}^{\infty}x_{n}$ converges.\\
Let $z_n=yx_n$ for all $n\in \mathbb{N}$. 
Then $FS(\langle{z_n}\rangle_{n=1}^{\infty})\subseteq  A\cap (0, \epsilon)$ and $\sum_{n=1}^{\infty}z_{n}$ 
be a convergent series in $S$. Therefore $p\cdot q\in \varGamma_0(S)$ which completes the proof.\\

\textbf{Theorem 2.7}: Let $S$ be a dense subsemigroup of $((0, \infty), +)$ such that $S\cap(0,1)$
is a subsemigroup of $((0, 1), \cdot)$ and assume that for each $y\in S\cap(0,1)$ and each $x\in S$, $x/y$ and
$yx\in S$. Let $r\in \mathbb{N}$ and let $S\cap(0,1) = \bigcup_{i=i}^{r}B_{i}$. Then there is some 
$i\in \{1,2,.....,r\}$ such that $B_{i}$ is an $IP$-set near zero and $B_{i}$ is central in $(S\cap(0,1),\cdot)$.\\

\textbf{Proof}: By Lemma 2.6, $\varGamma_0(S)$ is a left ideal of $\beta(S\cap(0,1),\cdot)$. Pick a minimal
left ideal $L$ of $\beta(S\cap(0,1),\cdot)$ with $L\subseteq \varGamma_0(S)$ and pick $p=p\cdot p\in L$.
Note that $p\in K(\beta(S\cap(0,1),\cdot))$ with $p=p\cdot p$ so all of its member are central in 
($S\cap(0,1),\cdot))$ Pick $i\in \{1,2,.....,r\}$ such that $B_{i}\in p$. Then $B_{i}$ is central in $(S\cap(0,1),\cdot)$
and, since $p\in \varGamma_0(S)$,  $ B_{i}$ is an $IP$-set near zero.\\

\textbf{Definition 2.8}: Let $S$ be a dense subsemigroup of $((0, \infty), +)$. Define $M_0(S)=\{p\in \beta S: $
if $A\in p$ then $A$ is a central set near zero\}.\\

\textbf{Lemma 2.9}: Let $S$ be a dense subsemigroup of $((0, \infty), +)$ such that $S\cap(0,1)$
is a subsemigroup of $((0, 1), \cdot)$ and assume that for each $y\in S\cap(0,1)$ and each $x\in S$, $x/y$ and
$yx\in S$. Then $M_0(S)$ is a left ideal of $\beta(S\cap(0,1),\cdot)$.\\

\textbf{Proof}: See the proof \cite[Theorem 5.6]{Hindman2}.\\

\textbf{Theorem 2.10}: Let $S$ be a dense subsemigroup of $((0, \infty), +)$ such that $S\cap(0,1)$
is a subsemigroup of $((0, 1), \cdot)$ and assume that for each $y\in S\cap(0,1)$ and each $x\in S$, $x/y$ and
$yx\in S$. Let $r\in \mathbb{N}$ and let $S\cap(0,1) = \bigcup_{i=i}^{r}B_{i}$. Then there is some 
$i\in \{1,2,.....,r\}$ such that $B_{i}$ is central near zero and $B_{i}$ is central in $(S\cap(0,1),\cdot)$.\\

\textbf{Proof}: See the proof \cite[Theorem 5.6]{Hindman2}.\\

\textbf{Theorem 2.11}: Let $S$ be a dense subsemigroup of $((0, \infty), +)$ such that $S\cap(0,1)$
is a subsemigroup of $((0, 1), \cdot)$ and assume that for each $y\in S\cap(0,1)$ and each $x\in S$, $x/y$ and
$yx\in S$. If $A$ is syndetic in $(S\cap (0, 1), \cdot)$ then $A$ is central near zero.\\

\textbf{Proof}: Since $A$ is syndetic in $(S\cap (0, 1), \cdot)$, there exists 
$G \in \mathcal{P}_{f}((0, 1)\cap S)$ such that $S\cap(0,1)=\bigcup_ {t\in G}t^{-1}A$.
Now take an idempotent $p$ in $K(0^+(S))$. Choose $t\in G$ such that $t^{-1}A\in p$.
Thus $t^{-1}A$ is central near zero. So by  \cite[Lemma 4.8]{Hindman2}, $A$ is central near zero.\\

\textbf{Theorem 2.12}: Let $S$ be a dense subsemigroup of $((0, \infty), +)$ such that $S\cap(0,1)$
is a subsemigroup of $((0, 1), \cdot)$ and assume that for each $y\in S\cap(0,1)$ and each $x\in S$, $x/y$ and
$yx\in S$. If $A$ is piecewise syndetic in $(S\cap (0, 1), \cdot)$ then for each $\epsilon>0$, $l\in \mathbb{N}$
there exist $a, d\in S$ such that $\{a, a+d, .....,a+(l-1)d\}\subseteq A\cap (0,\epsilon)$.\\

\textbf{Proof}: Le $I=\{p\in \beta S:$ if  $A\in p$ then for each $l\in \mathbb{N}$,
there exist $a, d\in S$ such that $\{a, a+d, .....,a+(l-1)d\}\subseteq A$\}.
By  \cite[Theorem 4.11]{Hindman2}, $E(K(0^+(S))\subseteq I$ (where for a semigroup $S$, $E(S)=\{x\in S: x$ is an idempotent\}). Let $I_0=I\cap 0^+(S)$ . 
Then clearly $I_0\neq \O{}$.
By theorem 1.8, $0^+(S)$ is a two-sided ideal of $\beta(S\cap(0,1),\cdot)$.
To show that $I_0$ is a two-sided ideal of $\beta(S\cap(0,1), \cdot)$, it is enough to prove that $I$ is a two-sided 
 ideal of $\beta(S\cap(0,1), \cdot)$. To this end let $p\in I$ and $q\in \beta(S\cap(0,1), \cdot)$ .
Suppose $A\in q\cdot p$ then $\{x\in S: x^{-1}A\in p\}\in q$. Choose $x\in S\cap(0,1)$ such that
$x^{-1}A\in p$. Pick $a, d\in S$ such that $\{a, a+d, .....,a+(l-1)d\}\subseteq x^{-1}A$. Then
$\{ax, ax+dx, .....,ax+(l-1)dx\}\subseteq A$. Thus $q\cdot p\in I$. Also let $A\in p\cdot q$.
Then $B=\{x\in S: x^{-1}A\in q\}\in p$. Pick $a, d\in S$ such that $\{a, a+d, .....,a+(l-1)d\}\subseteq B$.
Then $\bigcap_{i=0}^{l-1}(a+id)^{-1}A\in q$. Take $s\in S\cap(0,1)$ such that $s\in \bigcap_{i=0}^{l-1}(a+id)^{-1}A$.
Then $\{as, as+ds, .....,as+(l-1)ds\}\subseteq A$.
Thus $p\cdot q\in I$. Therefore $I$ is a two-sided ideal of $\beta(S\cap(0,1), \cdot)$. \\
Let $A$ be a  piecewise syndetic set  in $(S\cap (0, 1), \cdot)$. 
Now by Theorem 1.2 $K(\beta (S\cap(0,1), \cdot))\cap cl A\neq\O{}$.
Since $I_0$ is a two-sided ideal of $\beta (S\cap(0,1), \cdot))$, $K(\beta (S\cap(0,1), \cdot))\subseteq I_0$ 
and hence $clA\cap I_0\neq\O{}$.
Now Choose $p\in I_0$ such that $A\in p$. Therefore for each $\epsilon>0$, $l\in \mathbb{N}$
there exist $a, d\in S$ such that $\{a, a+d, .....,a+(l-1)d\}\subseteq A\cap (0,\epsilon)$.\\

    We can define $IP_0$-set near zero for a dense subsemigroup of $((0, \infty), +)$ as 
    is defined on $(\mathbb{N}, +)$\\
    
\textbf{Definition 2.13}: Let $S$ be a dense subsemigroup of $((0, \infty), +)$ and $A\subseteq S$.
Then $A$ is siad to be an $IP_0$-set near zero if for each $m\in \mathbb{N}$ and $\epsilon>0$ there exists a 
finite sequence $\langle{y_n}\rangle_{n=1}^{m}$ of positive reals  such that
$FS(\langle{y_n}\rangle_{n=1}^{m})\subseteq A\cap (0,\epsilon)$.\\

\textbf{Theorem 2.14}: Let $S$ be a dense subsemigroup of $((0, \infty), +)$ such that $S\cap(0,1)$
is a subsemigroup of $((0, 1), \cdot)$ and assume that for each $y\in S\cap(0,1)$ and each $x\in S$, $x/y$ and
$yx\in S$. If $A$ is central in $(S\cap (0, 1), \cdot)$ then $A$ is an $IP_0$-set near zero.\\

\textbf{Proof}: Let $J_0=\{p\in \beta S: $ if $A\in p$ then $A$ is an $IP_0$-set near zero\}.\\
Notice that $\varGamma_0(S)\subseteq J_0$ and therefore $ J_0\neq\O{}$. We now claim that 
$J_0$ is a two-sided ideal of $\beta (S\cap(0,1), \cdot)$. To this end let $p\in J_0$ 
and $q\in \beta (S\cap(0,1), \cdot)$. Suppose $A\in q\cdot p$. Then $B=\{x\in S: x^{-1}A\in p\}\in q$.
Choose $x\in S\cap(0,1)$ such that $x^{-1}A\in p$. Then for each $m\in \mathbb{N}$ and $\epsilon>0$ there exists a 
finite sequence $\langle{y_n}\rangle_{n=1}^{m}$ of positive reals such that
$FS(\langle{y_n}\rangle_{n=1}^{m})\subseteq x^{-1}A\cap (0,\epsilon)$.\\
Thus $FS(\langle{xy_n}\rangle_{n=1}^{m})\subseteq A\cap (0,\epsilon)$ and therefore $q\cdot p\in J_0$.
Also let $A\in p\cdot q$. Then $B=\{x\in S: x^{-1}A\in q\}\in p$.
So for each $m\in \mathbb{N}$ and $\epsilon>0$ there exists a 
finite sequence $\langle{y_n}\rangle_{n=1}^{m}$ of positive reals such that
$F=FS(\langle{y_n}\rangle_{n=1}^{m})\subseteq B$.
Then $\bigcap_{y\in F}y^{-1}A\in q$. Choose $x\in \bigcap_{y\in F}y^{-1}A\cap (0, 1)$.
Then $FS(\langle{xy_n}\rangle_{n=1}^{m})\subseteq A\cap (0,\epsilon)$.
Therefore $p\cdot q\in J_0$. Thus $J_0$ is a two-sided ideal of $\beta (S\cap(0,1), \cdot)$.
Let $A$ be central in $(S\cap(0,1), \cdot)$. Then there exists an idempotent $p\in K(\beta (S\cap(0,1), \cdot))$
such that $A\in p$. Also since $J_0$ is a two-sided ideal of $\beta (S\cap(0,1), \cdot)$, $p\in J_0$.
Therefore $A$ is an $IP_0$-set near zero.\\

\section{Monochromatic Solution To $\sum_{t=1}^nx_{t}=\prod_{t=1}^ny_t$ Near Zero}
In Ref. \cite{Hindman1}, author generalized the affirmative answer to the question of P. Csikv\'{a}ri, A. S\'{a}rk\"{o}zy, and K. Gyarmati
regarding the monochromatic solution of $a+b=cd $, whenever the set of naturals $\mathbb N$ is finitely colored. \\

 \textbf{Theorem 3.1}: Let $r\in \mathbb{N}$ and let $\mathbb{N}=\bigcup_{i=i}^{r}A_{i}$.There exists
 $i\in \{1,2,.....,r\}$
 such that for each $m\in  \mathbb{N}$, \\
 (1) there exists an increasing sequence $\langle{y_n}\rangle_{n=1}^{\infty}$ with distinct finite products 
 such that $FP(\langle{y_n}\rangle_{n=1}^{\infty})\subseteq A_{i}$ and whenever $F \in \mathcal{P}_{f}(\mathbb{N})$, 
 there exists a strongly increasing sequence $\langle{x_n}\rangle_{n=1}^{m}$ such that $FS(\langle{x_n}\rangle_{n=1}^{m})\subseteq A_{i}$
 and $\sum_{n=1}^mx_{n}=\prod_{n\in F}y_n$ and\\
 
 (2) there exists a strongly increasing sequence $\langle{x_n}\rangle_{n=1}^{\infty}$ 
 such that\\ $FS(\langle{x_n}\rangle_{n=1}^{\infty})\subseteq A_{i}$ and whenever $F \in \mathcal{P}_{f}(\mathbb{N})$,
 there exists an increasing sequence $\langle{y_n}\rangle_{n=1}^{m}$ with distinct finite products 
 such that $FP(\langle{y_n}\rangle_{n=1}^{m})\subseteq A_{i}$ and $\prod_{n=1}^{m}y_n =\sum_{n\in F}x_{n}$.\\
 
 \textbf{Proof}: See the proof \cite[Theorem 5]{Hindman1}.\\
 
 In this section we extend this result for a dense subsemigroup of $((0, \infty), +)$ to near zero.\\
 To establish the main result we need the following lemma.\\ 
 \textbf{Lemma 3.2}:  Let $S$ be a dense subsemigroup of $((0, \infty), +)$ such that $S\cap(0,1)$
is a subsemigroup of $((0, 1), \cdot)$ and assume that for each $y\in S\cap(0,1)$ and each $x\in S$, $x/y$ and
$yx\in S$. Let $\langle{w_t}\rangle_{t=1}^{\infty}$ be a sequence in $S\cap(0,1)$ such that 
$\sum_{t=1}^{\infty}w_{t}$ converges. Then there exist 
 sequences $\langle{x_t}\rangle_{t=1}^{\infty}$ and $\langle{y_t}\rangle_{t=1}^{\infty}$ such that 
 $\langle{x_t}\rangle_{t=1}^{\infty}$ is 
 strictly decreasing (and therefore it is distinct), $\langle{y_t}\rangle_{t=1}^{\infty}$ is decreasing 
 and has distinct finite products,
 $FS(\langle{x_t}\rangle_{t=1}^{\infty})\subseteq FS(\langle{w_t}\rangle_{t=1}^{\infty})$ and
$FP\langle{y_t}\rangle_{t=1}^{\infty}\subseteq FP(\langle{w_t}\rangle_{t=1}^{\infty})$

\textbf{Proof}: We construct the sequence $\langle{x_t}\rangle_{t=1}^{\infty}$ inductively. Since $\sum_{t=1}^{\infty}w_{t}$ converges,
$\langle{w_t}\rangle_{t=1}^{\infty}$ converges to $0$. Hence we can construct a strictly decreasing subsequence of 
$\langle{w_t}\rangle_{t=1}^{\infty}$ in the following way. Let $x_1=w_1$. Then clearly $(0, w_1/2)$ contains infinitely many points of the sequence
$\langle{w_t}\rangle_{t=2}^{\infty}$. Let us take $x_2=w_l$ such that $w_l\in (0, w_1/2)$. Then clearly $x_1>x_2$. Inductively,
let $\langle{x_t}\rangle_{t=1}^{k}$ be a strictly decreasing sequence such that for each $t$, $x_t=w_{p_t}$ and $x_t\in (0, x_{t-1}/2)$ for some $p_t \in \mathbb{N}$
and $x_t>x_{t+1}$ for $t\in \{1,2,.....,k-1\}$. Clearly $FS(\langle{x_t}\rangle_{t=1}^{k})\subseteq FS(\langle{w_t}\rangle_{t=1}^{\infty})$.
Now $(0, x_k/2)$ contains infinitely many points of the sequence $\langle{w_t}\rangle_{t=1}^{\infty}$ hence we can choose 
min$\{t\in \mathbb{N}: w_t\in (0, x_k/2)\}=q$. Let $x_{k+1}=w_q$.\\

 By similar arguments, we construct the sequence $\langle{y_t}\rangle_{t=1}^{\infty}$ inductively. Let
 $\langle{y_k}\rangle_{k=1}^{n}$ be a strictly decreasing subsequence of $\langle{w_t}\rangle_{t=1}^{\infty}$ such that $y_k=w_{t_k}$
 with $t_1<t_2.....<t_n$. Let $E=FP(\langle{y_k}\rangle_{k=1}^{n})$ also let $\mu=$ min $E\cup \{u^{-1}v: u,v\in E\}$.
 Now $(0,\mu)$ contains infinitely many points of $\langle{w_t}\rangle_{t={t_n}+1}^{\infty}$. Let
 min$\{t\in \mathbb{N}:t\geq {t_n}+1, w_t\in (0, \mu)\}=r$. Let $y_{n+1}=w_r$. Therefore,
 $FP\langle{y_t}\rangle_{t=1}^{\infty}\subseteq FP(\langle{w_t}\rangle_{t=1}^{\infty})$. This completes the proof.\\ 
 
 \textbf{Theorem 3.3}: Let $S$ be a dense subsemigroup of $((0, \infty), +)$ such that $S\cap(0,1)$
is a subsemigroup of $((0, 1), \cdot)$ and assume that for each $y\in S\cap(0,1)$ and each $x\in S$, $x/y$ and
$yx\in S$. Let $r\in \mathbb{N}$, $\epsilon>0$ and let $S\cap (0,1)=\bigcup_{i=i}^{r}A_{i}$. 
There exists $i\in \{1,2,.....,r\}$ such that for each $m\in  \mathbb{N}$,\\
 (1) there exists a decreasing sequence $\langle{y_n}\rangle_{n=1}^{\infty}$ with distinct finite products 
 such that $FP(\langle{y_n}\rangle_{n=1}^{\infty})\subseteq A_{i}\cap (0, \epsilon)$ and 
 whenever $F \in \mathcal{P}_{f}(\mathbb{N})$, 
 there exists a strictly  decreasing sequence 
 $\langle{x_n}\rangle_{n=1}^{m}$ such that $FS(\langle{x_n}\rangle_{n=1}^{m})\subseteq A_{i}\cap (0, \epsilon)$
 and $\sum_{n=1}^mx_{n}=\prod_{n\in F}y_n$ and\\
 (2) there exists a strictly decreasing sequence $\langle{x_n}\rangle_{n=1}^{\infty}$ 
 such that\\ $FS(\langle{x_n}\rangle_{n=1}^{\infty})\subseteq A_{i}\cap (0, \epsilon)$ and whenever
 $F \in \mathcal{P}_{f}(\mathbb{N})$,
 there exists a decreasing sequence $\langle{y_n}\rangle_{n=1}^{m}$ with distinct finite products 
 such that $FP(\langle{y_n}\rangle_{n=1}^{m})\subseteq A_{i}\cap (0, \epsilon)$ and 
 $\prod_{n=1}^{m}y_n =\sum_{n\in F}x_{n}$.\\
 
 \textbf{Proof}: Pick $p\in \beta (S\cap(0,1), \cdot)$ such that, for every $A\in p$ there exist
 sequences $\langle{x_n}\rangle_{n=1}^{\infty}$ and $\langle{y_n}\rangle_{n=1}^{\infty}$
 with $FS(\langle{x_n}\rangle_{n=1}^{\infty})\subseteq A\cap (0, \epsilon)$ and 
 $FP(\langle{y_n}\rangle_{n=1}^{\infty})\subseteq A\cap (0, \epsilon)$ (In section 2 we have proved the existence of such a $p$).
 Pick $i\in \{1,2,.....,r\}$ such that $A_{i}\cap (0, \epsilon)\in p$.\\
 Let $m\in \mathbb{N}$ be given. Let $B_0=\{z\in A_{i}\cap (0, \epsilon):$ there exists a strictly
 decreasing sequence $\langle{x_n}\rangle_{n=1}^{m}$ such that
 $FS(\langle{x_n}\rangle_{n=1}^{m})\subseteq A_{i}\cap (0, \epsilon)$ and $z=\sum_{n=1}^mx_{n}$\}.
 Let $C_0= \{z\in A_{i}\cap (0, \epsilon):$ there exists a
 decreasing sequence $\langle{y_n}\rangle_{n=1}^{m}$ with distinct finite products such that
 $FP(\langle{y_n}\rangle_{n=1}^{m})\subseteq A_{i}\cap (0, \epsilon)$ and $z==\prod_{n=1}^my_{n}$\}\\
 We claim that $B_0\in p$. If possible, let $B_0\notin p$, in which case \\
 $(A_{i}\cap (0, \epsilon))\backslash B_0\in p$. Pick a sequence $\langle{x_n}\rangle_{n=1}^{\infty}$
 with \\
 $FS(\langle{x_n}\rangle_{n=1}^{\infty})\subseteq (A_{i}\cap (0, \epsilon))\backslash B_0$. By Lemma 3.2
 we may assume that $\langle{x_n}\rangle_{n=1}^{\infty}$ is strictly decreasing. 
 But then $\sum_{n=1}^mx_{n}\in B_0$, a contradiction. Similarly $C_0\in p$.\\
 For conclusion (1) pick a decreasing sequence $\langle{y_n}\rangle_{n=1}^{\infty}$
 with distinct finite products such that $FP(\langle{y_n}\rangle_{n=1}^{\infty})\subseteq B_0$.
 For conclusion (2) pick a strictly decreasing sequence $\langle{x_n}\rangle_{n=1}^{\infty}$
 with distinct finite products such that $FS(\langle{x_n}\rangle_{n=1}^{\infty})\subseteq C_0$.\\

 \textbf{Corollary 3.4}: Let $S$ be a dense subsemigroup of $((0, \infty), +)$ such that $S\cap(0,1)$
is a subsemigroup of $((0, 1), \cdot)$ and assume that for each $y\in S\cap(0,1)$ and each $x\in S$, $x/y$ and
$yx\in S$. Let $r\in \mathbb{N}$ and let $\mathbb{N}=\bigcup_{i=i}^{r}A_{i}$. There exists $i\in \{1,2,.....,r\}$
and $a$, $b$, $c$, and $d$ in $S$ such that $\{a, b, c, d\}\subseteq A_{i}\cap (0, \epsilon)$
with $a+b= cd$ and $a\neq b$.\\
\textbf{Proof}: It is the special case of the above theorem when $n=2$


\begin{thebibliography}{3}

\bibitem{Bergelson2} V. Bergelson, Ultrafilters, IP sets, dynamics and combinatorial number theory;
Contemporary Mathematics Volume 530, 2010.


\bibitem{Bergelson1} V. Bergelson and D. Glasscock, Interplay Between Notions Of Additive And 
Multiplicative Largeneess, arXiv: 1610.09771v1 [math.CO] 31 Oct 2016.

\bibitem{Bergelson3} V. Bergelson, N. Hindman, On {\bf IP*}-sets and central sets, Combinatorica 
{\bf 14} (1994), 269-277.

\bibitem{csikvari} P. Csikv\'{a}ri, A. S\'{a}rk\"{o}zy, and K. Gyarmati, Density and Ramsey type 
results on algebraic equations with restricted  solution sets, Manuscript.

\bibitem{Hindman1} N. Hindman, Monochromatic sums equal to products in $\mathbb{N}$, Integers 
{\bf 11A}(2011), Article 10.
 
\bibitem{Hindman2} N. Hindman, I. Leader, The Semigroup of Ultrafilters Near 0;
     Semigroup Forum 59(1999), 33-55. 

\bibitem{Hindman3} N. Hindman, Amir Maleki and Dona Strauss, Central Sets and Their Combinatorial 
Characterization, Journal of Combinatorial Theory, (Series A) {\bf74} (1996),188-208.

\bibitem{Hindman4} Neil Hindman and Dona Strauss, Algebra in the Stone-$\breve{C}$ech 
compactification - theory and application, W.de Gruyter and Co.,Berlin, 1998.

\end{thebibliography}
\end{document}